\documentclass[11pt]{article}
\usepackage{authblk}

\usepackage{geometry}
\usepackage{amsmath,amssymb,mathrsfs}
\usepackage{color}
\usepackage{graphicx}
\usepackage{soul}
\usepackage{bbding}
\usepackage{lineno,hyperref}
\usepackage{multirow, booktabs, stmaryrd}
\SetSymbolFont{stmry}{bold}{U}{stmry}{m}{n}
\usepackage{bm}
\usepackage{subcaption}

\newcommand\btau{\boldsymbol{\tau}}

\usepackage{bm}

\def\bx{\mathbf{x}}
\def\bn{\mathbf{n}}
\def\bu{\mathbf{u}}

\def\bg{\mathbf{g}}
\def\bF{\mathbf{F}}
\def\bX{\mathbf{X}}
\def\grad{\nabla}
\def\div{\nabla\cdot}
\def\laplace{\Delta}

\newcommand{\beq}{\begin{equation}}
\newcommand{\eeq}{\end{equation}}
\newcommand{\beqs}{\begin{eqnarray}}
\newcommand{\eeqs}{\end{eqnarray}}
\newcommand{\beqsn}{\begin{eqnarray*}}
\newcommand{\eeqsn}{\end{eqnarray*}}
\newcommand{\bary}{\begin{array}}
\newcommand{\eary}{\end{array}}

\newcommand{\blue}[1]{\textcolor{black}{#1}}
\newcommand{\red}[1]{\textcolor{black}{#1}}

\newcommand*{\jb}[1]{\llbracket #1 \rrbracket}

\DeclareMathOperator*{\btheta}{\boldsymbol{\theta}}

\title{A hybrid neural-network and MAC scheme for Stokes interface problems}

\author[1,2]{Che-Chia Chang}
\author[1]{Chen-Yang Dai}
\author[3,4]{Wei-Fan Hu}
\author[1,4]{Te-Sheng Lin}
\author[1]{Ming-Chih Lai}

\affil[1]{Department of Applied Mathematics, National Yang Ming Chiao Tung University, Hsinchu 30010, Taiwan}
\affil[2]{Institute of Artificial Intelligence Innovation, National Yang Ming Chiao Tung University, Hsinchu 30010, Taiwan}
\affil[3]{Department of Mathematics, National Central University, Taoyuan 32001, Taiwan}
\affil[4]{National Center for Theoretical Sciences, National Taiwan University, Taipei 10617, Taiwan}

\begin{document}

\maketitle

\begin{abstract}
In this paper, we present a hybrid neural-network and MAC (Marker-And-Cell) scheme for solving Stokes equations with singular forces on an embedded interface in regular domains. As known, the solution variables (the pressure and velocity) exhibit non-smooth behaviors across the interface so extra discretization efforts must be paid near the interface in order to have small order of local truncation errors in finite difference schemes. The present hybrid approach avoids such additional difficulty. It combines the expressive power of neural networks with the convergence of finite difference schemes to ease the code implementation and to achieve good accuracy at the same time. The key idea is to decompose the solution into singular and regular parts. The neural network learning machinery incorporating the given jump conditions finds the singular part solution, while the standard MAC scheme is used to obtain the regular part solution with associated boundary conditions. The two- and three-dimensional numerical results show that the present hybrid method converges with second-order accuracy for the velocity  and first-order accuracy for the pressure,  and it is comparable with the traditional immersed interface method in literature.\\

\end{abstract}

\section{Introduction\label{sec1}}
In this paper, we consider $d$-dimensional ($d=2$ or $3$) Stokes equations in a regular domain $\Omega \subseteq \mathbb{R}^d$, in which an embedded interface $\Gamma$ with codimension $d-1$ (assumed to be smooth and closed) separates the domain into $\Omega^-$ and $\Omega^+$, respectively. Denoting the interface position by $\bX$, an interfacial force $\bF(\bX)$ defined only along the interface $\Gamma$ is exerted to the surrounding fluid and affects the flow behavior accordingly. Assuming the viscosity of both fluid subdomains are identical, one can write the single Stokes fluid system with an external body force $\bg$ as follows.
\begin{align}
\begin{split}
-\grad p(\bx) + \mu\laplace\bu(\bx) + \int_\Gamma \bF(\bX) \delta^d(\bx-\bX) \mbox{ d}\bX + \bg(\bx) = \mathbf{0}, \quad & \bx\in\Omega, \label{Stokes_IB} \\
\div\bu(\bx) = 0, \quad & \bx\in\Omega, \\
\bu (\bx) = \bu_b(\bx), \quad & \bx\in\partial\Omega,
\end{split}
\end{align}
where $\bu(\bx)$ and $p(\bx)$ are the velocity and the pressure, respectively; $\mu$ is the constant viscosity, and $\bu_b(\bx)$ is the velocity boundary condition. Notice that, the force term appeared in the first equation is singular and expressed in the Immersed Boundary (IB) formulation~\cite{P02} in which the integral involves a $d$-dimensional Dirac delta function $\delta^d$ over a $(d-1)$-dimensional surface resulting in one-dimensional delta function singularity. The above system~(\ref{Stokes_IB}) can be solved efficiently by the IB method~\cite{P02, TP09}. That is, the integral involving the delta function $\delta^d$ (line integral for $d = 2$ and surface integral for $d = 3$) can be regularized via a discrete delta function (a regularized form of the Dirac delta function) so the interfacial force $\bF$ can be spread into the fluid grid points near the interface. However, this singular force spreading process results in first-order accuracy for the velocity \cite{Mor08} and has $O(1)$ error for the pressure \cite{CFKL11}.

Due to the delta function singularity in Eq.~(\ref{Stokes_IB}), the pressure and velocity are no longer smooth across the interface so that the problem can be reformulated as the immersed interface formulation~\cite{LL97}
\begin{align}
    -\nabla p(\bx) + \mu \Delta \mathbf{u}(\bx) + \mathbf{g}(\bx) = \mathbf{0}, \quad & \bx \in \Omega^- \cup \Omega^+, \label{Stokes_4}\\
    \nabla \cdot \mathbf{u}(\bx) = 0, \quad & \bx \in \Omega^- \cup \Omega^+, \label{Stokes_5}\\
    \mathbf{u}(\mathbf{x}) = \mathbf{u}_b(\mathbf{x}), \quad  & \bx \in \partial \Omega, \label{Stokes_6}
\end{align}
where the pressure and velocity field are subjected to the following jump conditions (see the derivation in~\cite{LL01})
\begin{align}
    \llbracket p(\bX) \rrbracket = F_n(\bX), \quad &\bX \in \Gamma, \label{jump_p} \\
    \llbracket \bu(\bX) \rrbracket = \mathbf{0}, \quad \mu\left\llbracket \frac{\partial \mathbf{u}}{\partial \mathbf{n}}(\bX)\right\rrbracket = -(\bF(\bX) - F_n(\bX)\bn(\bX)), \quad  &\bX \in \Gamma. \label{jump_u}
\end{align}
Here, $F_n(\bX) = \bF(\bX)\cdot\bn(\bX)$ denotes the normal component of the interfacial force with $\bn(\bX)$ being the unit outward normal vector at $\bX\in\Gamma$. We use the double bracket \(\llbracket \cdot \rrbracket\) to denote the jump of a quantity  evaluated by the quantity from the \(\Omega^+\) side minus the one from the \(\Omega^-\) side. From the jump condition~(\ref{jump_p}), one can see that the pressure is discontinuous across $\Gamma$ when $F_n\neq0$. Also, from the jump condition~(\ref{jump_u}), the velocity is continuous across the interface while its normal derivative is discontinuous and determined by the tangential part of $\bF$. As a result, the velocity has a cusp behavior across the interface $\Gamma$.

As mentioned above, the pressure and velocity exhibit non-smooth behaviors across the interface, so extra discretization efforts must be paid near the interface in order to have small order of local truncation errors in finite difference schemes. A commonly used method is the Immersed Interface Method~\cite{LL97} (IIM) that incorporates the jump conditions into finite difference discretization on a uniform Cartesian grid and hereby achieves second-order convergence. Various versions of IIM can be found in~\cite{DZLYZ23, LIL12, LT08, WT22} and references therein. However, the implementation of jump incorporations could be tedious when dealing with complex interface geometry, especially in three-dimensional problems.

Until recent years, a new approach called Physics-Informed Neural Networks (PINNs)~\cite{RPK19}, was proposed to solve partial differential equations using machine learning methodology. The idea of PINN is to find a neural network approximate solution that minimizes the mean squared error loss consisting of the residuals of underlying differential equations along with boundary and initial conditions. Inspired by PINN, the authors have developed a series of structure-preserving neural network methods to solve elliptic interface problems \cite{HLL22, TLHL23} and Stokes interface problems \cite{TL23} that are able to capture the solution and derivative jump discontinuities sharply. These neural network solvers are mesh-free, so the correction of local truncation errors near the interface can be completely avoided. The results (in terms of accuracy) are also comparable with the ones obtained from traditional IIM mentioned above. \red{Other related works on the elliptic interface problem include the deep unfitted Nitsche method proposed by Guo et al.~\cite{GY22} and the convergence of PINNs method shown by Wu et al.~\cite{WZTL23}}.
Nevertheless, so far all the neural network solvers do not show a clean rate of convergence observed in traditional finite difference and finite element methods.
So in \cite{HLTL23}, the authors have developed a novel hybrid method that inherits both advantages of neural network approach and traditional finite difference  method  to solve Poisson interface problems. The entire computation only comprises a supervised learning task for function approximation and a fast direct solver for the Poisson equation, which can be easily and directly implemented regardless of the interface geometry. More importantly, this hybrid method shows second-order convergence in the numerical results. In this paper, we aim to extend this hybrid methodology for Poisson interface problem to the Stokes interface problem~(\ref{Stokes_4})-(\ref{Stokes_6}) with the jump conditions~(\ref{jump_p})-(\ref{jump_u}). Again, the whole computation of the proposed method simply consists of a supervised learning task and a regular Stokes solver. We should emphasize that we do not attempt to make a rigorous competition or comparison with existing traditional methods. Instead, we provide an alternative easy-to-implement way to solve Stokes interface problems by leveraging the strengths of machine learning approach and finite difference method.

The rest of the paper is organized as follows. In Section 2, we present the hybrid methodology with detailed implementations. Numerical results for Stokes interface problems are given in Sections 3, followed by some concluding remarks and future work in Section 4.

\section{Hybrid methodology for Stokes interface problems}

From the jump conditions (\ref{jump_p})-(\ref{jump_u}), it is clear to see that both the pressure and velocity of the Stokes equations are  piecewise-smooth due to the presence of the interface $\Gamma$. The key idea of the present hybrid approach is to decompose the solution variables into two distinct parts; namely, the singular and regular parts as
\begin{align}
    p(\bx)          & = p_r(\bx) + p_s(\bx),                   \label{p_decouple} \\
    \mathbf{u}(\bx) & = \mathbf{u}_r(\bx) + \mathbf{u}_s(\bx), \label{u_decouple}
\end{align}
where the subscript $r$ and $s$ denote the regular (smooth) and singular (non-smooth) component, respectively. More precisely, we hope that the regular part solutions $p_r$ and $\bu_r$ are fairly smooth over the entire fluid domain $\Omega$ so no jumps occur  across the interface. That is, we need to impose $\jb{p_r} = 0$ and $ \jb{\mathbf{u}_r} = \mu\jb{ \frac{\partial \mathbf{u}_r}{\partial \mathbf{n}}} = \mathbf{0}$  across the interface. On the other hand, the singular part solution variables $p_s$ and $\bu_s$ are responsible for having all discontinuous contributions; thus, from the jump conditions~(\ref{jump_p})-(\ref{jump_u}) and the decompositions~(\ref{p_decouple})-(\ref{u_decouple}), we can easily derive that
\begin{align}
    \llbracket p_s(\bX) \rrbracket = F_n(\bX), \quad &\bX \in \Gamma \label{jump_ps}, \\
    \llbracket \bu_s(\bX) \rrbracket = \bm{0}, \quad \mu\left\llbracket \frac{\partial \mathbf{u}_s}{\partial \mathbf{n}}(\bX)\right\rrbracket = -(\bF(\bX) - F_n(\bX)\bn(\bX))\equiv -\bF_\tau(\bX). \quad  &\bX \in \Gamma, \label{jump_us}
\end{align}
For succinctness, we use the term $\bF_\tau$ to represent the tangential force. In the following subsections, we shall explain how to obtain the singular and regular solutions in a sequential manner using our proposed hybrid methodology.

\subsection{Singular part solution}\label{subsec:singular}

The singular part solution is found by a supervised learning for neural network function approximations as follow. As mentioned earlier, both $p_s$ and the derivatives of $\bu_s$ should be discontinuous across the interface. To this end, we formally define the singular part solution in a piecewise manner as
\begin{align}
    p_s(\bx) =
    \begin{cases}
        \mathcal{P}(\bx), & \text{ if } \bx\in\Omega^{-}, \\
        0,                         & \text{ if } \bx\in\Omega^{+}, \\
    \end{cases} \quad \mbox{and} \quad
    \mathbf{u}_s(\bx) =
    \begin{cases}
        \mathbf{\mathcal{U}}(\bx), & \text{ if } \bx\in\Omega^{-}, \\
        \mathbf{0},                         & \text{ if } \bx\in\Omega^{+}, \\
    \end{cases}
    \label{decomposition}
\end{align}
where $\mathcal{P}$ and $\bm{\mathcal{U}}$ are neural network functions to be found. Plugging the above expression into the jump conditions~(\ref{jump_ps})-(\ref{jump_us}), we immediately derive that the unknown function $\mathcal{P}$ and $\bm{\mathcal{U}}$ must fulfill the following constraints along the interface:
\begin{align}
\mathcal{P}(\bX) = -F_n(\bX), \quad &\bX \in \Gamma, \label{P_net} \\
    \bm{\mathcal{U}}(\bX)=\bm{0}, \quad
    \mu \frac{\partial \bm{\mathcal{U}}(\bX)}{\partial \mathbf{n}} = \mathbf{F}_\tau(\bX), \quad
    -\nabla \mathcal{P}(\bX) + \mu \Delta \bm{\mathcal{U}}(\bX) - \jb{\bg(\bX)} = \bm{0}, \quad & \bX \in \Gamma. \label{U_net}
\end{align}
Note that, the third constraint of Eq.~(\ref{U_net}) is requested to have sufficient smoothness for the regular part solution (see Eq.~(\ref{Stokes_regular_1})). We shall examine this additional condition later in the next subsection.

In this paper, we utilize the standard supervised learning methodology to find the above unknown functions $\mathcal{P}$ and $\bm{\mathcal{U}}$ representing by the corresponding neural networks. We must point out that, it is not necessary to train a network model consisting of a $(\mathcal{P}, \bm{\mathcal{U}})$ pair to fulfill Eqs.~(\ref{P_net}) and (\ref{U_net}) simultaneously, although it seems that both of them are coupled together due to the third constraint of Eq.~(\ref{U_net}). In fact, such $\mathcal{P}$ and $\bm{\mathcal{U}}$ are not unique in the sense that there exists infinitely many functions defined in the whole domain $\Omega$ that satisfy the restrictions~(\ref{P_net}) and (\ref{U_net}) along $\Gamma$. This observation provides us some kind  of freedom to determine $\mathcal{P}$ and $\bm{\mathcal{U}}$ separately. Thus, we can first learn $\mathcal{P}$ to fulfill Eq.~(\ref{P_net}), and then use the obtained $\mathcal{P}$ to learn $\bm{\mathcal{U}}$ through Eq.~(\ref{U_net}) accordingly.

Let us illustrate the full training procedure in detail as follows. Firstly, we determine $\mathcal{P}$ to fulfill the restriction~(\ref{P_net}) via a Mean Squared Error (MSE) loss model. That is, given a dataset with $M$ training points $\{\bX^i\in\Gamma\}_{i=1}^{M}$, the corresponding loss function is defined by the MSE of the condition~(\ref{P_net}) as
\begin{align}
    \text{Loss}_\mathcal{P}(\btheta) = \frac{1}{M} \sum_{i=1}^{M} (\mathcal{P}(\bX^i; \btheta) + F_n(\bX^i))^2,
    \label{Loss_p}
\end{align}
where $\btheta$ collects all trainable parameters (weights and biases) of $\mathcal{P}$. Once $\mathcal{P}$ is trained, we learn the vector function $\bm{\mathcal{U}}$ that satisfies the constraints in Eq.~(\ref{U_net}). (Notice that $\bm{\mathcal{U}} =  (\mathcal{U}_1, \mathcal{U}_2)$ for $d = 2$ and $\bm{\mathcal{U}} =  (\mathcal{U}_1, \mathcal{U}_2, \mathcal{U}_3)$ for $d = 3$ case.) \blue{One could build a multi-output network to represent $\bm{\mathcal{U}}$; however, training the network is challenging because more terms are used in the whole loss function.} Instead, here, we construct independent single-output networks to learn each component of $\bm{\mathcal{U}}$. This step allows us to obtain all $\mathcal{U}_j$'s using parallel computations since all of them are decoupled (see Eq.~(\ref{U_net})). Using the same training data points, we train each $\mathcal{U}_j$ again by the MSE of Eq.~(\ref{U_net}) as
\begin{align}
\begin{split}
    \text{Loss}_{\mathcal{U}_j}(\boldsymbol{\theta}_j) = \frac{1}{M} \sum_{i=1}^{M} & \biggl[\biggl({\mathcal{U}_j}(\bX^i;\boldsymbol{\theta}_j)\biggr)^2
  + \biggl(\mu \frac{\partial \mathcal{U}_j}{\partial \bn}(\bX^i;\boldsymbol{\theta}_j)- F_{\tau_j}(\bX^i)\biggr)^2              \\
  & + \biggl(- \frac{\partial \mathcal{P}}{\partial x_j}(\bX^{i}; \bm{\theta}) + \mu \Delta \mathcal{U}_j(\bX^i; \bm{\theta}_j) - \jb{g_j(\bX^i)}\biggr)^2 \biggr],
\end{split}
\label{Loss_u}
\end{align}
where $\btheta_j$ denotes the trainable parameters of $\mathcal{U}_j$; $F_{\tau_j}$ and $g_j$ represent the $j$-th element of $\bF_\tau$ and $\bg$, respectively. One should note that, although we train the above network functions using only training data points along $\Gamma$, the trained functions are defined in the whole domain thanks to the expressive power of the neural networks. The spatial partial derivatives of the network functions appearing in the above loss models can be directly evaluated  using automatic differentiation. We employ the Levenberg-Marquardt method~\cite{Marquardt63}, a particularly designed efficient optimizer for nonlinear least-squares problems, to minimize the loss models~(\ref{Loss_p}) and (\ref{Loss_u}). With both $\mathcal{P}$ and $\bm{\mathcal{U}}$ (thus the singular part solution $p_s$ and $\bu_s$) in hand, we are ready to find the regular part solution  in the following subsection.

\subsection{Regular part solution}\label{subsec:regular}

Once we obtained the singular part solution via the neural network learning machinery, our next step is to solve the regular parts for $p_r$ and $\bu_r$. Substituting the decompositions~(\ref{p_decouple}) and (\ref{u_decouple}) into the Stokes equations~(\ref{Stokes_4})-(\ref{Stokes_6}), we immediately find that the regular part solution must satisfy the following Stokes-like system
\begin{align}
    -\nabla p_r(\bx) + \mu \Delta \bu_r(\bx) = \nabla p_s(\bx) - \mu \Delta \bu_s(\bx) - \mathbf{g}(\bx), \quad \quad & \bx\in\Omega, \label{Stokes_regular_1} \\
    \nabla \cdot \bu_r(\bx) = -\nabla \cdot \bu_s(\bx), \quad \quad & \bx\in\Omega, \label{Stokes_regular_2}\\
    \bu_r(\bx) = \bu_b(\bx), \quad \quad  & \bx \in \partial \Omega, \label{Stokes_regular_3}
\end{align}
where the boundary condition for $\bu_r$ (\ref{Stokes_regular_3}) is straightforwardly derived using Eq.~(\ref{Stokes_6}) and the definition of $\bu_s$ in Eq.~(\ref{decomposition}). Based on the above system, it is apparent to see that the solution pair $(\bu_r,p_r)$ satisfies the Stokes-like equations but the regular velocity is not divergence-free. \blue{Since we do not require the singular part solution to be divergence-free, the regular part solution must be compromised so that the resulting velocity is divergence-free.} In order to obtain smooth pressure and velocity  for the regular part, the right-hand side of  Eq.~(\ref{Stokes_regular_1}) must exhibit zero jump discontinuity across the interface which can be easily seen from the third constraint in Eq.~(\ref{U_net}). Meanwhile, the right-hand side of Eq.~(\ref{Stokes_regular_2}) also has zero jump across the interface since $\div \bm{\mathcal{U}} =0$ on $\Gamma$. To see this, one can write the divergence of $\bm{\mathcal{U}}$  as
\[
\div \bm{\mathcal{U}} = \grad_s \cdot \bm{\mathcal{U}} + \frac{\partial \bm{\mathcal{U}}}{\partial \mathbf{n}} \cdot \bn, \]
where the notation $\grad_s$ represents the surface gradient operator involving the derivatives along the surface $\Gamma$. From the first constraint of $\bm{\mathcal{U}}=0$ on $\Gamma$ in Eq.~(\ref{U_net}), we can conclude that all surface derivatives are zero leading to $\grad_s \cdot \bm{\mathcal{U}}=0$. And, from the second constraint in Eq.~(\ref{U_net}), the normal derivative of $\bm{\mathcal{U}}$ has only tangential components leading to $\frac{\partial \bm{\mathcal{U}}}{\partial \mathbf{n}} \cdot \bn=0$ as well. Therefore, we have $\div \bm{\mathcal{U}} = 0$ on the interface.

Now, the regular system~(\ref{Stokes_regular_1})-(\ref{Stokes_regular_3}) can be discretized by the well-known MAC scheme~\cite{HW65} and solved by the Uzawa-type algorithm~\cite{TAU92}. More precisely, the resultant linear system is a saddle point problem as
\begin{align}
    \begin{bmatrix}
        L      & G      \\
        G^\top & \bm{0}
    \end{bmatrix}
    \begin{bmatrix}
        \mathbf{u}_r \\ p_r
    \end{bmatrix}
    = \begin{bmatrix}
        \nabla p_s - \mu\Delta \mathbf{u}_{s} - \mathbf{g} \\ -\nabla \cdot \mathbf{u}_s
    \end{bmatrix}
    +
    \begin{bmatrix}
        \mathbf{f}(\mathbf{u}_b) \\ h(\mathbf{u}_b)
    \end{bmatrix}
    ,
    \label{linear_system}
\end{align}
where $L$ denotes the usual five-point discrete Laplacian operator $\mu \Delta_h$ applying to $\bu_r$, $G$ is the central difference gradient operator $-\nabla_h$ to $p_r$, and $G^\top$ is the central difference divergence operator $\nabla_h \cdot$ to $\bu_r$. Note that, those difference operators are all performed on the MAC grid. \blue{In the right-hand side of the matrix equation~(\ref{linear_system}), the first term involving the computation of derivatives of the singular part solution at those MAC grid (only inside the interface needed due to the definition~(\ref{decomposition})) can be done by automatic differentiation~\cite{GW08}. For instance, the term $\Delta \bu_s$ can be computed as follows.  Once $\bm{\mathcal{U}}$ is learned, the singular part solution $\bu_s$ is obtained. By the definition of $\bu_s$, we thus have $\Delta \bu_s(\bx)= \bm{0}, \bx \in\Omega^{+}$, and $\Delta \bu_s(\bx)= \Delta \bm{\mathcal{U}}(\bx), \bx \in\Omega^{-}$, where $\Delta \bm{\mathcal{U}}$ can be calculated directly by automatic differentiation.} The second term is associated with boundary condition of $\bu_r$ induced by the discretization operators $L$ and $G^\top$. Here, we use the linear interpolation to approximate the velocity $\bu_r$ at the ghost points  near the boundary. Notice that, the MAC scheme has the advantage that no boundary condition for the regular pressure $p_r$ is needed.

Rewriting the above linear system by
\begin{align*}
    \begin{bmatrix}
        L      & G      \\
        G^\top & \bm{0}
    \end{bmatrix}
    \begin{bmatrix}
        \mathbf{u}_r \\
        p_r
    \end{bmatrix}
    =
    \begin{bmatrix}
        b_1 \\
        b_2
    \end{bmatrix},
\end{align*}
the solution to the linear system can be found by applying the Schur complement technique as given in \cite{LHL12}. Using the block LU decomposition, one can further rewrite the system as
\begin{align}
    \begin{bmatrix}
        L      & \bm{0}          \\
        G^\top & -G^\top L^{-1}G
    \end{bmatrix}
    \begin{bmatrix}
        I      & L^{-1} G \\
        \bm{0} & I
    \end{bmatrix}
    \begin{bmatrix}
        \mathbf{u}_r \\
        p_r
    \end{bmatrix}
    =
    \begin{bmatrix}
        b_1 \\
        b_2
    \end{bmatrix},
    \label{eq:linear_system_2}
\end{align}
where \(I\) is the identity matrix. Then, we can solve the above linear system \eqref{eq:linear_system_2} via the following fractional steps by first introducing the intermediate velocity $ \bu^*=\bu_r+(L^{-1} G)p_r$.
\begin{enumerate}
    \item Solve $L\mathbf{u}^* = b_1$ to obtain the intermediate velocity field $\bu^*$. Notice that, this comprises a fully decoupled linear system of Poisson equation for each component of $\bu^*$, so we need to apply two ($d = 2$) or three ($d = 3$) fast Poisson solvers that can be efficiently done using the public software package such as Fishpack~\cite{ASS80} or Fast Fourier Transform (FFT).

    \item Solve $(-G^\top L^{-1} G) p_r = \tilde{b}_2 - G^\top\bu^*$ using the conjugate gradient method. In each iteration step, a matrix-vector multiplication $(-G^\top L^{-1} G)\psi$ for some vector $\psi$ is required; the inversion of $L^{-1}$ again can be performed efficiently by the fast Poisson solvers mentioned in Step 1. We must emphasize that, since the kernel of the underlying matrix $-G^\top L^{-1} G$ is obviously spanned by the constant vector $\bm{1} = [1,1,\cdots,1]^\top$, here we use the modified right-hand side vector $\tilde{b}_2 =  b_2 - \bm{1}^\top b_2/\|\bm{1}\|^2$ instead of the original  $b_2$. That is, we project $b_2$ into the solution space to guarantee the existence of the solution $p_r$. The iteration stops whenever the system residual is less than $10^{-12}$.

    \item Compute the regular part velocity $\bu_r = \bu^* - (L^{-1} G)p_r$. Again, this involves the usage of the fast Poisson solvers.
\end{enumerate}
As a conclusion for the computational complexity, the overall cost in the above three steps can be counted in terms of the number of fast Poisson solvers being applied.

\subsection{Summary of the hybrid method}

Let us summarize the step by step procedure of the proposed hybrid method as follows.
\begin{enumerate}
    \item Use neural network learning machinery to construct $\mathcal{P}$ and $\bm{\mathcal{U}}$, and form the singular part solution $p_s$ and $\bu_s$.
    \item Apply the MAC scheme to discretize the Stokes-like system and use the Uzawa-type algorithm to  find the regular part solution for the velocity  $\bu_r$ and pressure $p_r$.
    \item Add the singular and regular parts of solution to obtain the desired velocity $\bu$ and pressure $p$ to the Stokes interface problem.
\end{enumerate}

Let us conclude this section by remarking several features of the present method in the following. The proposed method builds upon the hybrid method introduced in \cite{HLTL23} for Poisson interface problem and extends to solve the Stokes interface problem. In contrast to traditional sharp interface methods such as in \cite{LL97}, the present method avoids the effort of careful constructing additional correction terms near the interface grid points to achieve desired accuracy.  We utilize a supervised learning technique for the singular part solution and traditional finite difference scheme for the regular part solution.  Once the network solution \(\mathcal{P}\) and \(\bm{\mathcal{U}}\) have been trained, they can be utilized to solve the problem for various grid sizes. This allows adapting the method to different resolutions flexibly. Another notable advantage is the ability to leverage well-established fast Poisson solvers, which facilitates the efficient solution of linear systems arising in this method. This efficiency is crucial for practical implementation and enables handling larger and more complex problems. Furthermore, it is straightforward to extend the proposed method to multiple-interface problems.
We shall also point out that the numerical error of the present method is induced by the neural network approximations (optimization and approximation errors) and finite difference discretizations (truncation error). 
Despite that, we can observe the numerical convergence behavior in our results shown in the next section which is hardly seen in complete neural network approaches.

\section{Numerical results}

In this section, we perform a series of numerical tests to validate the efficiency and convergence of the proposed method. For all examples, we fix the constant viscosity by $\mu = 1$. 
 We use the sigmoid function as the activation function for the neural networks to approximate \(\mathcal{P}\) and \(\bm{\mathcal{U}}\), and stop the network training once the loss value reaches the threshold $10^{-10}$. All the problems considered here are defined on regular domains with a uniform layout of $N^d$ grid points.
The numerical experiments were conducted on a PC equipped with an Intel Core i7-13700 CPU and 128GB RAM which gives sufficient computational resources for all the tests.

\paragraph{\textbf{Example 1.}}
The first example, referred from the 2D IIM Stokes solver in~\cite{C01}, aims to validate the convergence property of the present hybrid method. The domain is set as  $\Omega = [-2,2]^2$ with a unit circular interface, $\Gamma = \{\bX(\theta) = (\cos \theta,\sin \theta)|\theta\in[0,2\pi)\}$ immersed in the domain. The interfacial force consists of both tangential and normal components as $\bF(\bX(\theta)) = 2\sin(3\theta)\boldsymbol{\tau}(\theta) - \cos^3(\theta) \mathbf{n}(\theta)$ in which the tangent vector can be calculated by $\btau(\theta) = \bX'(\theta)=(-\sin \theta, \cos \theta)$ and the normal vector thus is $\bn(\theta) = (\cos \theta,\sin \theta)$.

The exact velocity $\bu = (u_1, u_2)$ written in polar coordinates is chosen by
\begin{align*}
     & u_1(r, \theta) =
    \begin{cases}
        \frac{1}{8}r^2\cos(2\theta) + \frac{1}{16}r^4\cos(4\theta) - \frac{1}{4}r^4 \cos(2\theta)             & \text{ if } r<1,      \\
        -\frac{1}{8} r^{-2} \cos(2\theta) + \frac{5}{16} r^{-4}\cos(4\theta) - \frac{1}{4}r^{-2}\cos(4\theta) & \text{ if } r \geq 1,
    \end{cases} \\
     & u_2(r, \theta)    =
    \begin{cases}
        - \frac{1}{8} r^2 \sin(2\theta) + \frac{1}{16} r^4 \sin(4\theta) + \frac{1}{4}r^4 \sin(2\theta)         & \text{ if } r<1,     \\
        \frac{1}{8} r^{-2} \sin(2\theta) + \frac{5}{16} r^{-4} \sin(4\theta) - \frac{1}{4} r^{-2} \sin(4\theta) & \text{ if } r\geq 1,
    \end{cases}
\end{align*}
while the exact pressure written in Cartesian coordinates ($x = r\cos\theta, y = r\sin\theta$) is given as
\[
p(x,y) =
\begin{cases}
x^3 + \cos(\pi x) \cos(\pi y) & \text{ if } r<1, \\
\cos(\pi x)\cos(\pi y)        & \text{ if } r \geq 1.
\end{cases}
\]
With the above exact velocity and pressure, we can obtain the external force $\bg = (g_1, g_2)$ in Cartesian coordinates accordingly as
\begin{align*}
     & g_1(x,y) =
    \begin{cases}
        - \pi \sin(\pi x)\cos(\pi y)+ 6x^2 - 3y^2                                & \text{ if } r<1, \\
        - \pi \sin(\pi x)\cos(\pi y) -\frac{3(x^4 - 6x^2y^2 + y^4)}{(x^2+y^2)^4} & \text{ if } r\geq 1, \\
    \end{cases} \\
     & g_2(x,y) =
    \begin{cases}
        - \pi \sin(\pi x)\cos(\pi y)-  6xy                                    & \text{ if } r<1, \\
        - \pi \sin(\pi x)\cos(\pi y)- \frac{12 (x^3 y - xy^3)}{(x^2+y^2)^{4}} & \text{ if } r\geq 1.
    \end{cases}
\end{align*}

\blue{As discussed earlier,  we first learn the neural network functions $\mathcal{P}$ and $\bm{\mathcal{U}}$ using the losses of Eq.~(\ref{Loss_p}) and (\ref{Loss_u}), respectively. Then we use the learned functions to obtain the singular part solution for all various grid sizes. Theoretically, one might want to choose the number of training points on the interface $M$ to be linearly scaled with the grid number $N$ accordingly so that mesh resolutions on the interface and the domain can be consistent.  However, for the sake of efficiency, we simply choose the number $M$ to have a successful training (the loss value reaches the threshold $10^{-10}$). Here, we 
set the number $M=400$ based on the grid size $N=256$  so that mesh resolutions on the interface and the domain are almost the same. We train three independent shallow (one-hidden-layer) neural networks (one for $\mathcal{P}$ and the other two for $\bm{\mathcal{U}}=(\mathcal{U}_1,\mathcal{U}_2)$) with $50$ neurons in the hidden layer and $M=400$ randomly sampled training points along the interface $\Gamma$. The training results are shown in Table~\ref{tab:ex1_nn}. Within a few hundred epochs, the training losses decrease to below $10^{-10}$, and the training time take less than three seconds.}

\begin{table}[!ht]
	\centering
    \begin{tabular}{lccc}
        \toprule
        {}    & training time (s) & final training loss & epochs \\
        \midrule
        $\mathcal{P}$   & 0.13         & 8.949e-11           & 54     \\
        $\mathcal{U}_1$ & 2.63         & 9.878e-11           & 328    \\
        $\mathcal{U}_2$ & 2.94         & 3.920e-11           & 353    \\
        \bottomrule
    \end{tabular}
    \caption{Elapsed time for training the singular part solution in Example~1.}
    \label{tab:ex1_nn}
\end{table}

In Table~\ref{tab:ex1}, we report the $L^\infty$-errors as well as the corresponding rates of convergence of the numerical velocity 
$ \mathbf{u}= \mathbf{u}_r + \mathbf{u}_s =(u_1,u_2)$ and pressure $p=p_r+p_s$. 
 As seen clearly, all the fluid variables achieve approximately second-order convergence, and the magnitudes of error are comparable with the ones obtained by IIM~\cite{LT08}. \blue{Although this particular choice of number $M=400$ is based on the grid size used $N=256$, we further run the tests for the finer cases $N=512, 1024$, and the results keep the second-order convergence. It seems that the network solution for the singular part is sufficiently accurate so it does not pollute the overall accuracy. In other words, the rate of convergence is completely determined by the computation of the regular part solution in this example.}
 
\blue{Furthermore, we also check the divergence-free condition by computing the $L^\infty$-error of $\nabla \cdot \mathbf{u}$, denoting by $e_{\infty} (\nabla \cdot \mathbf{u}) = \| \nabla_h \cdot \mathbf{u}_r + \nabla \cdot \mathbf{u}_s\|_\infty$. Here, the regular part divergence is calculated by applying the usual finite difference scheme to $\mathbf{u}_r$ on MAC grid while the singular part divergence is calculated directly by applying the automatic differentiation to $\mathbf{u}_s$ at the same positions as the regular part. (In fact, we only need to compute the singular part divergence  $\nabla \cdot \mathbf{u}_s$ in the sub-domain $\Omega^-$ since it is zero in $\Omega^+$ by definition.) One can see that the $L^\infty$-error of $\nabla \cdot \mathbf{u}$ can be reached up to the order of magnitude $10^{-7}$ after the grid point $N=128$ is used. Thus, the present methodology preserves the divergence-free condition quite well.}

\begin{table}[h]
    \centering
    \begin{tabular}{lrrrrrrr}
    \toprule
    $N$   & $e_\infty(u_1)$ & rate  & $e_\infty(u_2)$ & rate & $e_\infty(p)$ &  rate & $e_\infty(\nabla \cdot \mathbf{u})$ \\
    \midrule
        $32$  & 2.714e-02       & -     & 2.062e-02       & -     & 9.804e-02     & -     &  3.742e-05  \\
        $64$  & 6.108e-03       & 2.15  & 4.685e-03       & 2.14  & 3.488e-02     & 1.49  & 1.209e-06 \\
        $128$ & 1.553e-03       & 1.98  & 1.256e-03       & 1.90  & 8.477e-03     & 2.04  & 3.580e-07\\
        $256$ & 3.578e-04       & 2.12  & 2.870e-04       & 2.13  & 2.418e-03     & 1.81  & 5.664e-07 \\
    \midrule    
        $512$ & 8.952e-05       & 2.00  & 7.571e-05       & 1.92  & 6.015e-04     & 2.01  & 5.866e-07  \\
        $1024$ & 2.116e-05       & 2.08  & 2.140e-05       & 1.82  & 1.547e-04     & 1.96  & 5.899e-07  \\
    \bottomrule
\end{tabular}
    \caption{\blue{Mesh refinement results of the 2D Stokes interface problem on $N^2$ grid points in Example~1.}}
    \label{tab:ex1}
\end{table}

\paragraph{\textbf{Example 2.}}

This example will showcase the reliability of our hybrid method in the scenario when  the exact solution is unavailable. In this test, an elliptical interface $\Gamma = \{\bX(\theta) = (0.5\cos\theta, 0.3\sin \theta)|\theta\in[0,2\pi)\}$ is immersed in the square domain $\Omega = [-1,1]^2$ while the interfacial force comprises a curvature normal force and tangential force as $\bF(\bX(\theta)) = 0.1\kappa(\theta) \mathbf{n}(\theta) - 0.1\bm{\tau}(\theta)$, where $\kappa(\theta)$ denotes the local curvature at $\bX(\theta)$. The fluid boundary condition $\bu_b$ is obtained by using the spectrally accurate boundary integral method~\cite{SL10} with free-space Stokeslets kernel (i.e., $\bu(\bx)\rightarrow\bm{0}$ as $\|\bx\|\rightarrow\infty$).

We follow the same setup for the network training procedure as in Example~1; namely, we train three independent shallow neural networks with $50$ neurons in the hidden layer and $400$ randomly sampled training points along the interface $\Gamma$. The training results are shown in Table~\ref{tab:ex2_nn}. In fewer than $1000$ epochs, the training losses decrease to below $10^{-10}$, and the training process takes less than eight seconds.

\begin{table}[!ht]
    \centering
    \begin{tabular}{lccc}
        \toprule
        {}    & training time(s) & final training loss & epochs \\
        \midrule
        $\mathcal{P}$   & 0.28         & 4.870e-11           & 150    \\
        $\mathcal{U}_1$ & 7.64         & 8.742e-11           & 878    \\
        $\mathcal{U}_2$ & 7.86         & 8.118e-11           & 906    \\

        \bottomrule
    \end{tabular}
    \caption{Elapsed time for training the singular part solution in Example~2.}
    \label{tab:ex2_nn}
\end{table}

To measure the convergence of the solutions, since the explicit solution is not available, we compute the error in a successive manner with $\|\phi_N - \phi_{2N}\|_\infty$, where $\phi_N$ denotes some fluid variable solution with the grid number $N$. It is worth mentioning that the interpolation on the MAC grid layout is not straightforward  since the solutions at different resolutions would not coincide at the same location. In the present method, we simply apply the linear interpolation to the regular part solution
and use neural network approximation for the singular part solution at target grid points.

The mesh refinement results are shown in Table~\ref{tab:ex2}. When increasing the grid number, all numerical solutions converge with roughly second-order accuracy. In addition, we compute the interfacial velocity $(u_1^\Gamma, u_2^\Gamma)$ at the discrete marker points $\{(0.5\cos(\theta_k),0.3\sin(\theta_k))| \theta_k = 2k\pi/64, k = 1, 2,\cdots, 64\}$. Again, this is done by the interpolation technique aforementioned. Referring the exact solution  obtained by the boundary integral method with spectral accuracy~\cite{SL10}, Table~\ref{tab:ex2_interface} indicates that the numerical interfacial velocity obtained by the present method also attains second-order convergence.

\begin{table}[!]
    \centering
    \begin{tabular}{lrrrrrr}
    \toprule
    $N$    & $e_\infty(u_1)$ & rate & $e_\infty(u_2)$ & rate & $e_\infty(p)$ & rate \\
    \midrule
        $32$  & 1.638e-03       & -     & 1.002e-03       & -     & 2.267e-02     & -     \\
        $64$  & 4.465e-04       & 1.87  & 2.579e-04       & 1.96  & 6.514e-03     & 1.80  \\
        $128$ & 1.071e-04       & 2.06  & 6.371e-05       & 2.02  & 1.635e-03     & 1.99  \\
        $256$ & 2.698e-05       & 1.99  & 1.604e-05       & 1.99  & 5.843e-04     & 1.48  \\
    \bottomrule
    \end{tabular}
    \caption{Mesh refinement results of the 2D Stokes interface problem on $N^2$ grid points in Example 2.}
    \label{tab:ex2}
\end{table}

\begin{table}[!]
    \centering
    \begin{tabular}{lrrrr}
    \toprule
    N    & $e_\infty(u_1^{\Gamma})$ & rate & $e_\infty(u_2^{\Gamma})$ & rate \\
    \midrule
        $32$  & 4.507e-04                & -     & 6.773e-04                & -     \\
        $64$  & 9.893e-05                & 2.19  & 1.352e-04                & 2.33  \\
        $128$ & 2.666e-05                & 1.89  & 2.117e-05                & 2.67  \\
        $256$ & 6.565e-06                & 2.02  & 4.518e-06                & 2.23  \\
    \bottomrule
    \end{tabular}
    \caption{Mesh refinement results of the interfacial velocity $(u_1^\Gamma,u_2^\Gamma)$ in Example 2. We use the results of boundary integral method~\cite{SL10} as the exact solution and compute the errors.}
    \label{tab:ex2_interface}
\end{table}

\begin{figure}[!]
    \centering
    \begin{subfigure}[b]{0.48\textwidth}
        \centering
        \includegraphics[width=\textwidth]{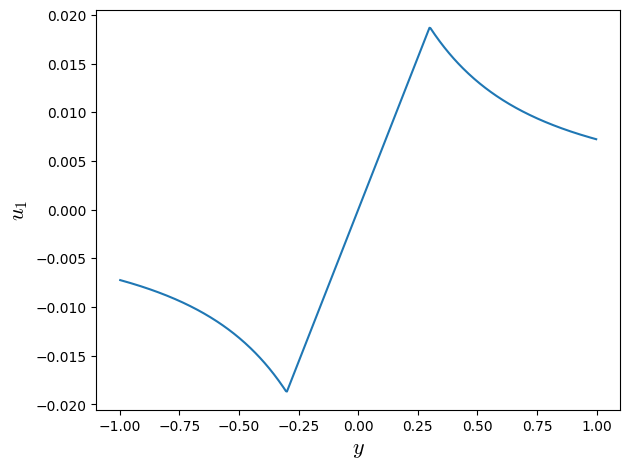}
        \caption{\(u_1\) along the grid line \(x=0\).}
    \end{subfigure}
    \begin{subfigure}[b]{0.48\textwidth}
        \centering
        \includegraphics[width=\textwidth]{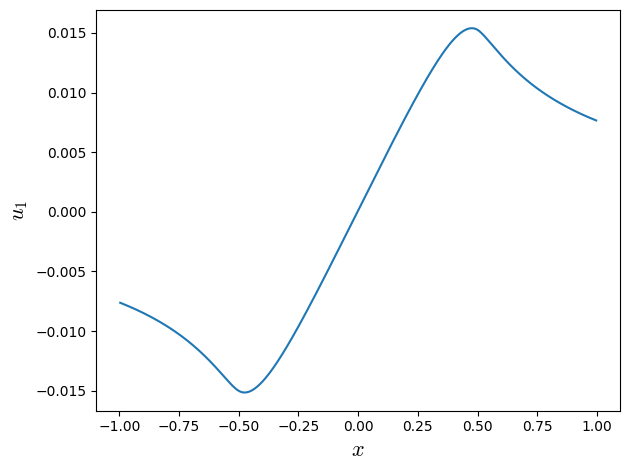}
        \caption{\(u_1\) along the grid line \(y=h/2\).}
    \end{subfigure}
    \begin{subfigure}[b]{0.48\textwidth}
        \centering
        \includegraphics[width=\textwidth]{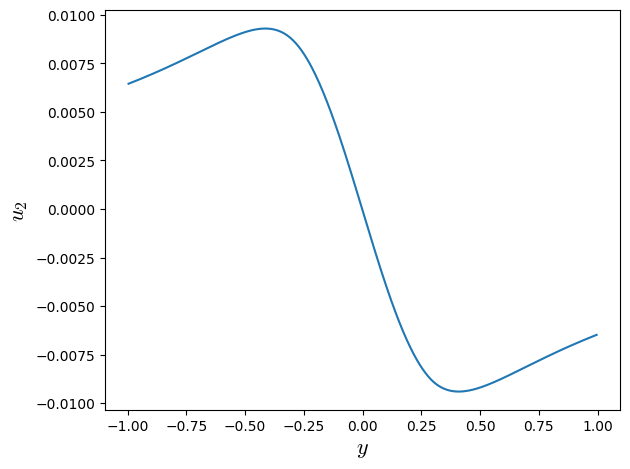}
        \caption{\(u_2\) along the grid line \(y=0\).}
    \end{subfigure}
    \begin{subfigure}[b]{0.48\textwidth}
        \centering
        \includegraphics[width=\textwidth]{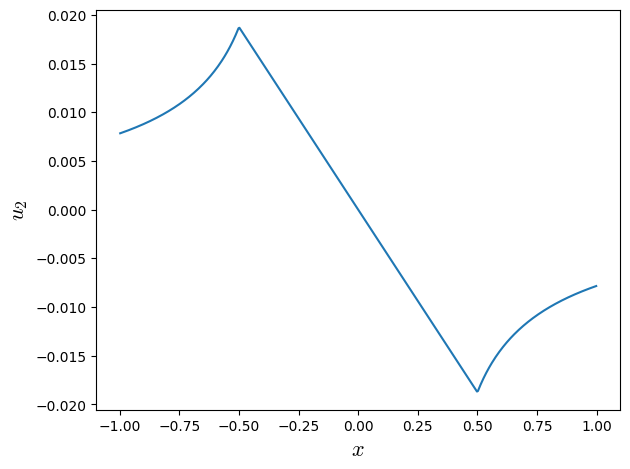}
        \caption{\(u_2\) along the grid line \(x=h/2\).}
    \end{subfigure}
    \begin{subfigure}[b]{0.48\textwidth}
        \centering
        \includegraphics[width=\textwidth]{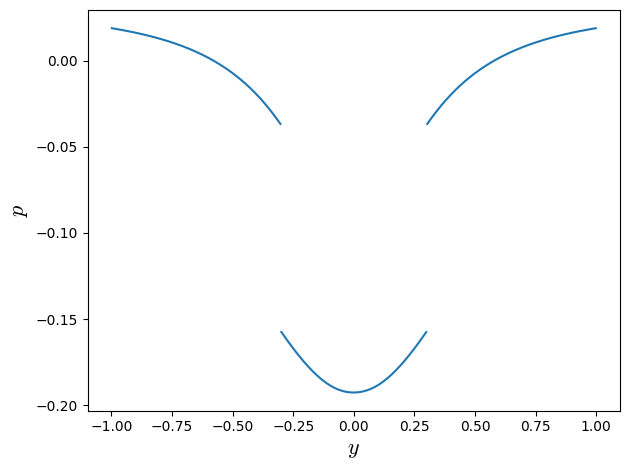}
        \caption{\(p\) along the grid line \(x=h/2\).}
    \end{subfigure}
    \begin{subfigure}[b]{0.48\textwidth}
        \centering
        \includegraphics[width=\textwidth]{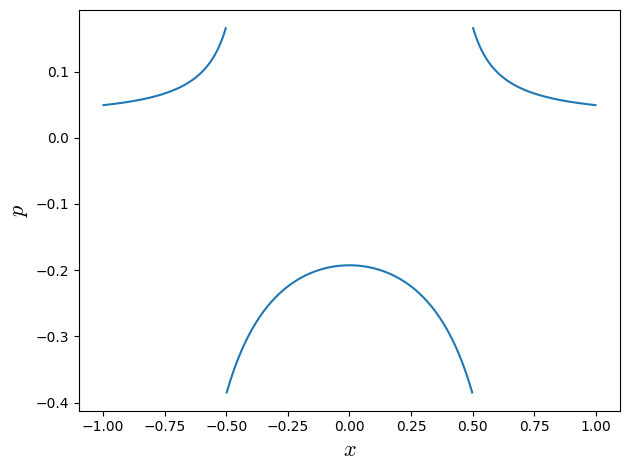}
        \caption{\(p\) along the grid line \(y=h/2\).}
    \end{subfigure}

    \caption{The cross-sectional plots of the velocity field and pressure along the grid lines in Example 2. The mesh width is $h=2/N$ with $N=512$.}
    \label{fig:cross_section}
\end{figure}

Lastly, to illustrate the ability of our method in capturing the derivative discontinuities for the velocity field and the jump discontinuity for the pressure, we show the cross-sectional views of the numerical velocity field and pressure in Fig.~\ref{fig:cross_section}. It can be clearly seen that the cusp (derivative discontinuity) for the velocity and jump for the pressure occur at the interface.

\paragraph{\textbf{Example 3.}}

Now, we proceed to solve the test example referred from~\cite{WT22} for three-dimensional Stokes interface problem in the cubic domain $\Omega = [-2,2]^3$. The interface \(\Gamma\) is simply described by the unit sphere $\Gamma = \{\bX = (x,y,z) \mid x^2+y^2+z^2 = 1\}$. Here, we choose the exact velocity $\bu = (u_1,u_2,u_3)$ as
\begin{align*}
     & u_1(x,y,z) =
    \begin{cases}
        \frac{1}{4}yz              & \text{ if } \mathbf{x} \in \Omega^-, \\
        \frac{1}{4}yz(x^2+y^2+z^2) & \text{ if } \mathbf{x} \in \Omega^+,
    \end{cases} \\
     & u_2(x,y,z) =
    \begin{cases}
        \frac{1}{4}xz              & \text{ if } \mathbf{x} \in \Omega^-, \\
        \frac{1}{4}xz(x^2+y^2+z^2) & \text{ if } \mathbf{x} \in \Omega^+,
    \end{cases} \\
     & u_3(x,y,z) =
    \begin{cases}
        -\frac{1}{2}xy(1-x^2-y^2) & \text{ if } \mathbf{x} \in \Omega^-, \\
        -\frac{1}{2}xyz^2        & \text{ if } \mathbf{x} \in \Omega^+,
    \end{cases}
\end{align*}
and the exact pressure $p$
\begin{align*}
    p(x,y,z) =
    \begin{cases}
        (-\frac{3}{4}x^3 + \frac{3}{8}x)yz & \text{ if } \mathbf{x} \in \Omega^-, \\
        0                                  & \text{ if } \mathbf{x} \in \Omega^+.
    \end{cases}
\end{align*}
One can accordingly obtain the external force field $\bg = (g_1, g_2, g_3)$ given by
\begin{align*}
     & g_1(x,y,z) =
    \begin{cases}
        (-\frac{9}{4}x^2 + \frac{3}{8})yz & \text{ if } \mathbf{x} \in \Omega^-, \\
        -\frac{7}{2}yz                    & \text{ if } \mathbf{x} \in \Omega^+,
    \end{cases} \\
     & g_2(x,y,z) =
    \begin{cases}
        (-\frac{3}{4}x^3 + \frac{3}{8}x)z & \text{ if } \mathbf{x} \in \Omega^-, \\
        -\frac{7}{2}xz                    & \text{ if } \mathbf{x} \in \Omega^+,
    \end{cases} \\
     & g_3(x,y,z) =
    \begin{cases}
        (-\frac{3}{4}x^3 - \frac{45}{8}x)y & \text{ if } \mathbf{x} \in \Omega^-, \\
        xy                                 & \text{ if } \mathbf{x} \in \Omega^+.
    \end{cases}
\end{align*}
The interfacial force can be derived from the exact velocity and pressure given above.  Also, the boundary conditions are given by the exact velocity.

We use a one-hidden-layer neural network with $50$ neurons  for the pressure and one-hidden-layer neural network with $100$ neurons for each component of the velocity to train the networks. The training data is generated by choosing $1000$ random points on the sphere. The training results are shown in Table~\ref{tab:ex3_nn}. Again, the pressure only takes less than one second to train, while each velocity component takes about $10-20$ seconds.

\begin{table}[!ht]
    \centering
    \begin{tabular}{lccc}
        \toprule
        {}    & training time(s) & final training loss & epochs \\
        \midrule
        $\mathcal{P}$   & 0.43         & 7.347e-11           & 181    \\
        $\mathcal{U}_1$ & 13.02        & 8.118e-11           & 803    \\
        $\mathcal{U}_2$ & 12.07        & 6.640e-11           & 764    \\
        $\mathcal{U}_3$ & 19.06        & 8.268e-11           & 1246   \\

        \bottomrule
    \end{tabular}
    \caption{Elapsed time for training the singular part solution in Example~3.}
    \label{tab:ex3_nn}
\end{table}

In Table~\ref{tab:ex3}, we show the \(L^\infty\)-errors and rates of convergence of the numerical solutions. The velocity shows clear second-order convergence, while the pressure is reduced to first-order only. Although not shown here, the largest numerical errors for the computed pressure occur at the corners of the cube. We believe that employing the quadratic interpolation to approximate the velocity $\bu_r$ at the ghost points near the boundary could  help reduce this error.

\begin{table}[!]
    \centering
    \begin{tabular}{lrrrrrrrr}
    \toprule
    $N$   & $e_\infty(u_1)$ & rate & $e_\infty(u_2)$ & rate & $e_\infty(u_3)$ & rate & $e_\infty(p)$ & rate \\
    \midrule
        $16$  & 3.626e-02       & -     & 3.626e-02       & -     & 2.258e-02       & -     & 1.283e+00     & -     \\
        $32$  & 1.041e-02       & 1.80  & 1.041e-02       & 1.80  & 6.534e-03       & 1.79  & 7.553e-01     & 0.76  \\
        $64$  & 2.768e-03       & 1.91  & 2.768e-03       & 1.91  & 1.773e-03       & 1.88  & 4.137e-01     & 0.87  \\
        $128$ & 7.123e-04       & 1.96  & 7.123e-04       & 1.96  & 4.653e-04       & 1.93  & 2.177e-01     & 0.93  \\
    \bottomrule
\end{tabular}
    \caption{Mesh refinement results of the 3D Stokes interface problem in Example 3.}
    \label{tab:ex3}
\end{table}

\paragraph{\textbf{Example 4.}}

In this example, we solve a three-dimensional Stokes interface problem with an oblate spheroid as the interface. The domain is a cube \(\Omega = [-1,1]^3\). The interface is given by \(\Gamma = \{\bX=(x,y,z) \mid \frac{x^2}{0.5^2} + \frac{y^2}{0.5^2} + \frac{z^2}{0.3^2} = 1\}\) that can be parameterized by \(\bX(\theta,\varphi) = (0.5\cos\varphi\cos\theta, 0.5\cos\varphi\sin\theta, 0.3\sin\varphi)\), where \(\theta \in [0,2\pi)\) and \(\varphi \in [-\frac{\pi}{2},\frac{\pi}{2}]\). Here, we simply take the interfacial force as the surface tension force (with the constant surface tension to be one) as  $\bm{F} =  - 2 H \mathbf{n}$, where  \(H\) is the mean curvature and \(\mathbf{n}\) is the unit normal vector. The no-slip boundary condition $\bu_b(\bx) = \mathbf{0}$ is applied along the boundary domain $\partial\Omega$ and the zero external force field \(\mathbf{g} = \bm{0}\) is assumed.

We use a one-hidden-layer neural network with $100$ neurons  for the pressure and two-hidden-layer neural network with $30$ neurons for each component of the velocity to train the networks. The training points are generated by choosing $1000$ random points on the interface. The training results are shown in Table~\ref{tab:ex6_nn}.

\begin{table}[!ht]
    \centering
    \begin{tabular}{lccc}
        \toprule
          & training time(s) & final training loss & epochs \\
        \midrule
        $\mathcal{P}$   & 10.097        & 1.850e-16           & 3000   \\
        $\mathcal{U}_1$ & 121.154       & 8.999e-13           & 2519   \\
        $\mathcal{U}_2$ & 142.336       & 2.346e-12           & 3000   \\
        $\mathcal{U}_3$ & 144.442       & 1.399e-12           & 3000   \\
        \bottomrule
    \end{tabular}
    \caption{Elapsed time for training the singular part solution in Example~4.}
    \label{tab:ex6_nn}
\end{table}

\begin{figure}[!]
    \centering
    \includegraphics[scale=0.3]{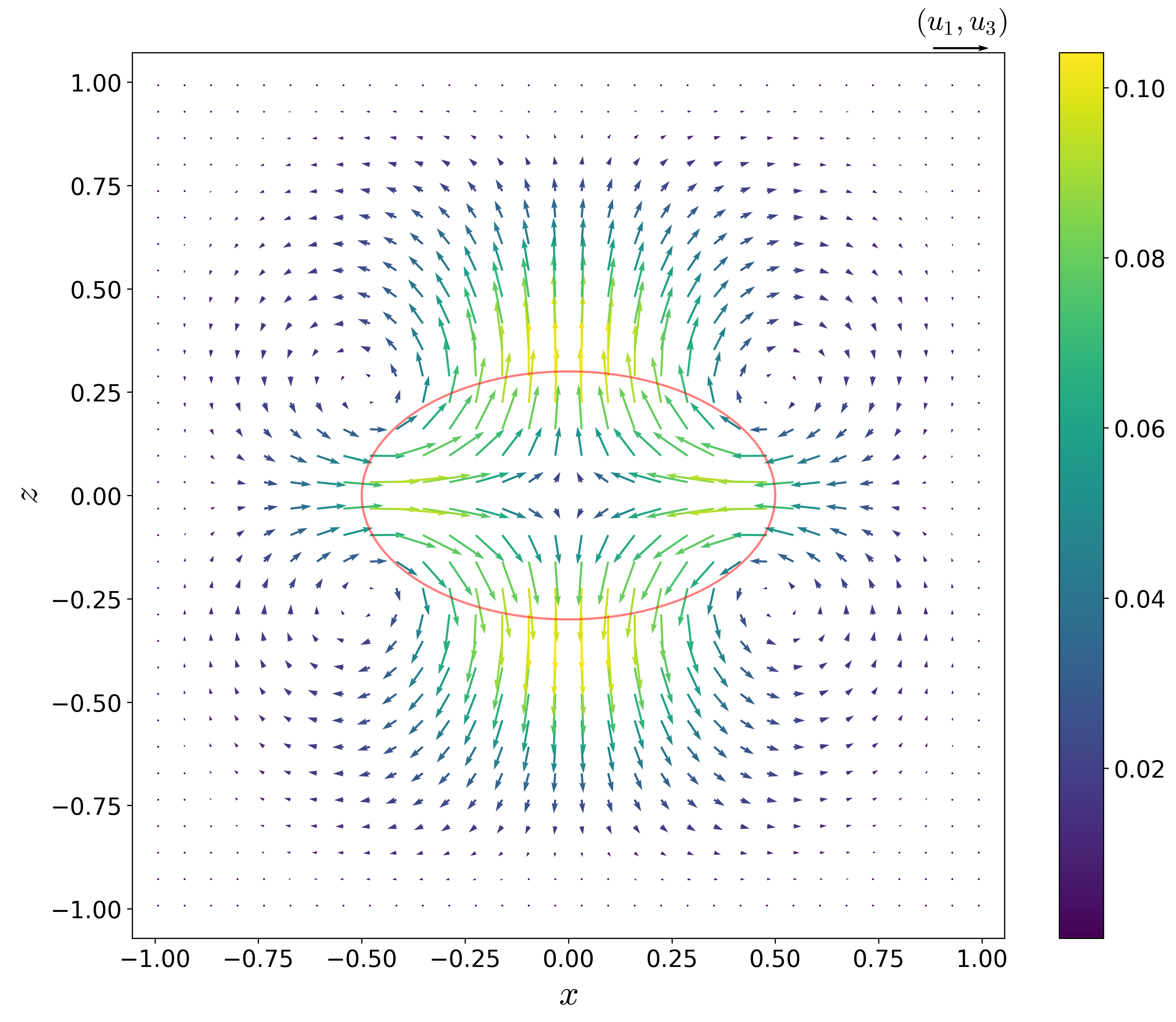}
    \caption{Quiver plot for $(u_1,u_3)$ along the plane $y=h/2$ in Example~4. The mesh width is $h=2/N$ with $N=128$. The color indicates the absolute magnitude of the velocity field.}
    \label{fig:ex6_quiver}
\end{figure}

As in Example 2, the explicit formulas for the solution are not available so we compute the errors in the same manner as there. In Table~\ref{tab:ex6}, we see that the numerical results show second-order convergence in the velocity, while first-order convergence in the pressure. We further depict the cross-sectional view of the flow quiver $(u_1,u_3)$ along the plane $y = 2/h$ in Fig.~\ref{fig:ex6_quiver}. One can see that, the flow tends to reduce the absolute magnitudes of the mean curvature along the interface and relaxes to a spherical shape. This instant flow tendency matches well with the simulations when the interface dynamics is considered. Meanwhile, since the fluid is incompressible, we can see that two vortex dipoles (or four counter-rotating vortices) occur near both sides. Therefore, our proposed hybrid method is indeed able to predict physically reasonable results.
\begin{table}[!ht]
    \centering
    \begin{tabular}{lrrrrrrrr}
        \toprule
    $N$   & $e_\infty(u_1)$ & rate & $e_\infty(u_2)$ & rate & $e_\infty(u_3)$ & rate & $e_\infty(p)$ & rate \\
        \midrule
    $32$  & 7.694e-03     & -  & 5.531e-03     & -  & 1.131e-02     & - & 1.057e+00     & - \\
    $64$  & 2.096e-03     & 1.88  & 1.418e-03     & 1.96  & 3.169e-03     & 1.84  & 5.731e-01     & 0.88  \\
    $128$ & 5.429e-04     & 1.95  & 4.326e-04     & 1.71  & 8.681e-04     & 1.87  & 2.945e-01     & 0.96  \\
        \bottomrule
    \end{tabular}
    \caption{Mesh refinement results of the 3D Stokes interface problem in Example 4.}
    \label{tab:ex6}
\end{table}

\section{Conclusion and future work}

In this paper, we present a hybrid neural-network and finite-difference method for solving Stokes equations with singular forces on an interface in regular domains. Overall speaking, this hybrid method takes both advantages of neural network and finite difference scheme. Our approach is based on decomposing the solution into two parts: the singular part captures the non-smooth solution behavior while the regular part represents the smooth solution. To find the singular part, we make use of neural network function approximators to learn the solution variables that satisfy certain constraints along the interface. Once the singular part solution is obtained, the regular part can be solved via a Stokes-like system using traditional MAC scheme. The overall computational costs can be calculated in terms of the number of fast Poisson solvers used in the iterations.  Through several numerical experiments, we demonstrate that the present hybrid method indeed obtains numerical results with traditional convergence property for both two- and three-dimensional Stokes interface problems. As a future work, we plan to apply the present method to solve time-dependent dynamic interface problems for flows in different applications, such as the vesicle hydrodynamics, electro-hydrodynamics, or phoretic systems.

\section*{Acknowledgements}

W.-F. Hu, T.-S. Lin, and M.-C. Lai acknowledge the supports by National Science and Technology Council, Taiwan, under the research grants 111-2115-M-008-009-MY3, 111-2628-M-A49-008-MY4, and 110-2115-M-A49-011-MY3, respectively. W.-F. Hu and T.-S. Lin also acknowledge the supports by National Center for Theoretical Sciences, Taiwan.


\end{document}